\documentclass{article}
\usepackage{amssymb}
\usepackage{amsfonts}
\usepackage{amsmath}

\newtheorem{theorem}{Theorem}

\newtheorem{corollary}[theorem]{Corollary}

\newtheorem{definition}[theorem]{Definition}

\newtheorem{remark}[theorem]{Remark}

\newenvironment{proof}[1][Proof]{\noindent\textbf{#1.} }{\ \rule{0.5em}{0.5em}}
\input{tcilatex}
\begin{document}

\title{From a dynamical system of the knee \\
to natural jet geometrical objects}
\author{Mircea Neagu and Mihaela Maria Marin}
\date{}
\maketitle

\begin{abstract}
In this paper we construct some natural geometrical objects on the $1$-jet
space $J^{1}(\mathbb{R},\mathbb{R}^{3})$, like a nonlinear connection, a
Cartan linear connection (together with its d-torsions and d-curvatures), a
jet "electromagnetic" d-field and its geometric "electromagnetic" Yang-Mills
energy, starting from a given dynamical system governing the
three-dimensional motion of the knee in the mathematical model introduced by
Grood and Suntay. The corresponding Yang-Mills energetic surfaces of
constant level (produced by this knee dynamical system) are studied.
\end{abstract}



\section{Short introduction}

The $1$-jet spaces are basic objects in the study of classical and quantum
field theories (see Olver \cite{Olver}). For such a reason, a lot of authors
(Asanov \cite{Asanov}, Saunders \cite{Saunders[8]} and many others) studied
the differential geometry of $1$-jet spaces. Using as a pattern the
Miron-Anastasiei's Lagrangian geometrical ideas (see \cite{Mir-An[4]}),
Balan and Neagu have recently developed the \textit{single-time Lagrange
geometry on }$1$\textit{-jet spaces} (see \cite{Balan-Neagu}), which is very
suitable for the geometrical study of the solutions of a given ODE system,
via the \textit{least squares variational method} initiated by Udri\c{s}te
(see \cite{Udr Geom Dyn[9]}).

In the present days, the biomechanics of the lower limb is an extremely
fruitful field of research. This research field is of interest for a lot of
biomechanists (e.g., see Barbu \cite{Barbu} and Ro\c{s}ca \cite{Rosca}). In
this direction, it is important to note that the differential geometrical
methods from $1$-jet spaces are very useful for studying biomechanics (see
Ivancevic \cite{Ivanc}). In such a perspective, this paper applies the jet
geometrical results from monograph \cite{Balan-Neagu} to the Euler's
dynamical equations that govern the three-dimensional motion of the tibia
with respect to the femur, in the biomechanical and mathematical model
introduced by Grood and Suntay \cite{Grood-Suntay}, and studied further by
many scholars (e.g., see Hefzy and Abdel-Rahman \cite{Hefzy}).

\section{Jet single-time Lagrange geometry produced by a non-linear
dynamical system}

Let us present now the main jet single-time Lagrangian geometrical results
that, in our opinion, may characterize a given non-linear dynamical system.
In this way, let us consider the jet fibre bundle of order one $J^{1}(%
\mathbb{R},\mathbb{R}^{n})\rightarrow \mathbb{R}\times \mathbb{R}^{n},$
where $n\geq 1,$ whose local coordinates $(t,x^{i},x_{1}^{i}),$ $i=\overline{%
1,n},$ obey the rules%
\begin{equation*}
\widetilde{t}=\widetilde{t}(t),\quad \widetilde{x}^{i}=\widetilde{x}%
^{i}(x^{j}),\quad \widetilde{x}_{1}^{i}=\frac{\partial \widetilde{x}^{i}}{%
\partial x^{j}}\frac{dt}{d\widetilde{t}}\cdot x_{1}^{j}.
\end{equation*}

Let $\mathfrak{X}=\left( \boldsymbol{X}_{(1)}^{(i)}(x^{k})\right) $ be an
arbitrary distinguished (d-) tensor field on the $1$-jet space $J^{1}(%
\mathbb{R},\mathbb{R}^{n})$, whose local components transform by the rules%
\begin{equation*}
\widetilde{\boldsymbol{X}}_{(1)}^{(i)}=\frac{\partial \widetilde{x}^{i}}{%
\partial x^{j}}\frac{dt}{d\widetilde{t}}\cdot \boldsymbol{X}_{(1)}^{(j)}.
\end{equation*}

The d-tensor field $\mathfrak{X}$ produces the jet first order ODE system (%
\textit{jet dynamical system})%
\begin{equation}
x_{1}^{i}=\boldsymbol{X}_{(1)}^{(i)}(x^{k}(t)),\quad \forall \text{ }i=%
\overline{1,n},  \label{DEs1}
\end{equation}%
where $c(t)=(x^{i}(t))$ is an unknown curve on $\mathbb{R}^{n}$ (i.e., a 
\textit{jet field line} of the d-tensor field $\mathfrak{X}$), and we have%
\begin{equation*}
x_{1}^{i}:=\dot{x}^{i}=\frac{dx^{i}}{dt},\quad \forall \text{ }i=\overline{%
1,n}.
\end{equation*}

Let us consider now the \textit{jet least squares Lagrangian function}
(attached to the dynamical system (\ref{DEs1}) and to Euclidian manifolds $(%
\mathbb{R},1)$ and $(\mathbb{R}^{n},\delta _{ij})$)%
\begin{equation*}
\mathbf{JLS}^{\text{ODEs}}:J^{1}(\mathbb{R},\mathbb{R}^{n})\rightarrow 
\mathbb{R}_{+},
\end{equation*}%
expressed by%
\begin{eqnarray}
\mathbf{JLS}^{\text{ODEs}}(x^{k},x_{1}^{k}) &=&\sum_{i,j=1}^{n}\delta _{ij}%
\left[ x_{1}^{i}-\boldsymbol{X}_{(1)}^{(i)}(x)\right] \left[ x_{1}^{j}-%
\boldsymbol{X}_{(1)}^{(j)}(x)\right] =  \notag \\
&=&\sum_{i=1}^{n}\left[ x_{1}^{i}-\boldsymbol{X}_{(1)}^{(i)}(x)\right] ^{2},
\label{JetLS}
\end{eqnarray}%
where $x=(x^{k})_{k=\overline{1,n}}.$ It is obvious that the global minimum
points of the \textit{jet least squares energy action}%
\begin{equation*}
\mathcal{E}^{\text{ODEs}}(c(t))=\int_{a}^{b}\mathbf{JLS}^{\text{ODEs}%
}(x^{k}(t),\dot{x}^{k}(t))dt,\quad t\in \lbrack a,b],
\end{equation*}%
are exactly the solutions of class $C^{2}$ of the jet dynamical system (\ref%
{DEs1}). In other words, any solution of class $C^{2}$ of the system (\ref%
{DEs1}) verifies the second order Euler-Lagrange equations (i.e., the 
\textit{jet geometric dynamics} associated to the ODE system (\ref{DEs1}))
produced by the jet least squares Lagrangian function (\ref{JetLS}): 
\begin{equation}
\frac{\partial \left[ \mathbf{JLS}^{\text{ODEs}}\right] }{\partial x^{i}}-%
\frac{d}{dt}\left( \frac{\partial \left[ \mathbf{JLS}^{\text{ODEs}}\right] }{%
\partial \dot{x}^{i}}\right) =0,\quad \forall \text{ }i=\overline{1,n}.
\label{E-L}
\end{equation}

\begin{remark}
Conversely, the preceding statement is not true. In other words, there exist
solutions of the Euler-Lagrange ODE system of second order \emph{(\ref{E-L})}%
, which are not global minimum points for the jet least squares energy
action $\mathcal{E}^{\text{\emph{ODEs}}}$, that is which are not solutions
for the initial jet dynamical system \emph{(\ref{DEs1})}.
\end{remark}

In such a context, we consider that we may regard the jet least squares
Lagrangian function $\mathbf{JLS}^{\text{ODEs}}$ as a natural geometrical
substitut on the $1$-jet space $J^{1}(\mathbb{R},\mathbb{R}^{n})$ for the
dynamical system (\ref{DEs1}). But, an entire single-time Lagrange geometry
on the $1$-jet space $J^{1}(\mathbb{R},\mathbb{R}^{n})$ (in the sense of
nonlinear connection, generalized Cartan linear connection, d-torsions,
d-curvatures, jet "electromagnetic" d-field and jet "electromagnetic"
Yang-Mills energy), geometry which is produced only by the jet least squares
Lagrangian function $\mathbf{JLS}^{\text{ODEs}}$ (via its Euler-Lagrange
equations (\ref{E-L})), is now completely done in the book \cite{Balan-Neagu}%
. For such a reason, we introduce the following concept:

\begin{definition}
Any kind of geometrical object on $J^{1}(\mathbb{R},\mathbb{R}^{n})$, which
is produced by the jet least squares Lagrangian function $\mathbf{JLS}^{%
\text{\emph{ODEs}}}$ (via its second order Euler-Lagrange equations \emph{(%
\ref{E-L})}) is called a \textbf{geometrical object produced by the jet
dynamical system} \emph{(\ref{DEs1})}.
\end{definition}

Let us consider the Jacobian matrix

\begin{equation*}
J\left( \mathfrak{X}\right) =\left( \dfrac{\partial \boldsymbol{X}%
_{(1)}^{(i)}}{\partial x^{j}}\right) _{i,j=\overline{1,n}}.
\end{equation*}%
In such a context, we give the following geometrical result (which is proved
in the book \cite{Balan-Neagu}):

\begin{theorem}
\label{MainTh}

\begin{enumerate}
\item[\emph{(1)}] The \textbf{canonical nonlinear connection} on $J^{1}(%
\mathbb{R},\mathbb{R}^{n})$\ \textbf{produced by the jet dynamical system }%
\emph{(\ref{DEs1})} has the local components 
\begin{equation*}
\Gamma ^{\text{\emph{ODEs}}}=\left( \boldsymbol{M}_{(1)1}^{(i)}=0,%
\boldsymbol{N}_{(1)j}^{(i)}\right) ,
\end{equation*}%
where%
\begin{eqnarray*}
\mathbf{N}_{(1)} &=&\left( \boldsymbol{N}_{(1)j}^{(i)}=-\frac{1}{2}\left[ 
\frac{\partial \boldsymbol{X}_{(1)}^{(i)}}{\partial x^{j}}-\frac{\partial 
\boldsymbol{X}_{(1)}^{(j)}}{\partial x^{i}}\right] \right) _{i,j=\overline{%
1,n}}= \\
&=&-\dfrac{1}{2}\left[ J\left( \mathfrak{X}\right) -\text{ }^{\emph{T}%
}J\left( \mathfrak{X}\right) \right] .
\end{eqnarray*}

\item[\emph{(2)}] All adapted components of the \textbf{canonical Cartan
linear connection} $\mathtt{C}\Gamma ^{\text{\emph{ODEs}}}$\ \textbf{%
produced by the jet dynamical system} \emph{(\ref{DEs1})} are zero.

\item[\emph{(3)}] The \textbf{torsion tensor} $\mathbf{T}^{\text{\emph{ODEs}}%
}$ of the canonical Cartan linear connection $\mathtt{C}\Gamma ^{\text{\emph{%
ODEs}}}$ \textbf{produced by the jet dynamical system} \emph{(\ref{DEs1})}
has the adapted components ($k=\overline{1,n}$)%
\begin{equation*}
\mathbf{T}_{(1)k}=\left( \boldsymbol{T}_{(1)jk}^{(i)}=-\frac{1}{2}\left[ 
\frac{\partial ^{2}\boldsymbol{X}_{(1)}^{(i)}}{\partial x^{k}\partial x^{j}}-%
\frac{\partial ^{2}\boldsymbol{X}_{(1)}^{(j)}}{\partial x^{k}\partial x^{i}}%
\right] \right) _{i,j=\overline{1,n}}=\dfrac{\partial \mathbf{N}_{(1)}}{%
\partial x^{k}}.
\end{equation*}

\item[\emph{(4)}] All adapted components of the \textbf{curvature tensor} $%
\mathbf{R}^{\text{\emph{ODEs}}}$ of the canonical Cartan linear connection $%
\mathtt{C}\Gamma ^{\text{\emph{ODEs}}}$ \textbf{produced by the jet
dynamical system} \emph{(\ref{DEs1})} cancel.

\item[\emph{(5)}] The \textbf{geometric "electromagnetic" distinguished }$2$%
\textbf{-form produced by the jet dynamical system} \emph{(\ref{DEs1})} has
the expression 
\begin{equation*}
\mathcal{F}^{\text{\emph{ODEs}}}=\boldsymbol{F}_{(i)j}^{(1)}\delta
x_{1}^{i}\wedge dx^{j},
\end{equation*}%
where 
\begin{equation*}
\delta x_{1}^{i}=dx_{1}^{i}+\boldsymbol{N}_{(1)k}^{(i)}dx^{k},\quad \text{ }%
\forall \text{ }i=\overline{1,n},
\end{equation*}%
and we have%
\begin{equation*}
\mathbf{F}^{(1)}=\left( \boldsymbol{F}_{(i)j}^{(1)}=\frac{1}{2}\left[ \frac{%
\partial \boldsymbol{X}_{(1)}^{(i)}}{\partial x^{j}}-\frac{\partial 
\boldsymbol{X}_{(1)}^{(j)}}{\partial x^{i}}\right] \right) _{i,j=\overline{%
1,n}}=-\mathbf{N}_{(1)}.
\end{equation*}

\item[\emph{(6)}] The \textbf{jet geometric "electromagnetic" Yang-Mills
energy produced by the jet dynamical system} \emph{(\ref{DEs1})} is given by
the formula 
\begin{equation*}
\mathbf{EYM}^{\text{\emph{ODEs}}}(x)=\sum_{i=1}^{n-1}\sum_{j=i+1}^{n}\left[ 
\boldsymbol{F}_{(i)j}^{(1)}\right] ^{2}=\dfrac{1}{2}\cdot \text{\emph{Trace}}%
\left[ \mathbf{F}^{(1)}\cdot \text{ }^{\emph{T}}\mathbf{F}^{(1)}\right] .
\end{equation*}
\end{enumerate}
\end{theorem}

\begin{remark}
The adapted components $\boldsymbol{F}_{(i)j}^{(1)}$ of the
"electromagnetic" (d-) $2$-form $\mathcal{F}^{\text{\emph{ODEs}}}$ produced
by the jet dynamical system \emph{(\ref{DEs1})} verify the following \textbf{%
geometrical Maxwell equations}:%
\begin{equation*}
\sum_{\{i,j,k\}}\boldsymbol{F}_{(i)j||k}^{(1)}=0,
\end{equation*}%
where $\sum_{\{i,j,k\}}$ represents a cyclic sum and we have%
\begin{equation*}
\boldsymbol{F}_{(i)j||k}^{(1)}=\frac{\partial \boldsymbol{F}_{(i)j}^{(1)}}{%
\partial x^{k}}.
\end{equation*}%
For more details in the jet "electromagnetic" topic, please consult the work 
\emph{\cite{Balan-Neagu}}.
\end{remark}

\begin{remark}
The jet geometric Yang-Mills energy $\mathbf{EYM}^{\text{\emph{ODEs}}}$
coincides with the norm of the skew-symmetric "electromagnetic" matrix $%
\mathbf{F}^{(1)}$ in the Lie algebra%
\begin{equation*}
o(n)=L(O(n))=\left\{ \left. A\in M_{n}(\mathbb{R})\text{ }\right\vert \text{ 
}A+\text{ }^{\emph{T}}A=0\right\} ,
\end{equation*}%
where $O(n)=\left\{ \left. A\in M_{n}(\mathbb{R})\text{ }\right\vert \text{ }%
A\cdot \text{ }^{\emph{T}}A=I_{n}\right\} $ is the corresponding Lie group
of the orthogonal matrices.
\end{remark}

\section{Three-dimensional mathematical model used for studying the
tibio-femoral dynamics}

The femur and tibia are modeled as two rigid bodies. Cartilage deformation
is assumed relatively small compared to joint motions and not to affect
relative motions and forces within the tibio-femoral joint. Furthermore,
friction forces will be neglected because of the extremely low coefficients
of friction of the articular surfaces. Hence, in this model, the resistance
to motion is essentially due to the ligamentous structures and the contact
forces. The menisci were not taken into consideration in the present model.
For more biomechanical details, see Hefzy and Abdel-Rahman \cite{Hefzy}.

The joint coordinate system, which was introduced by Grood and Suntay \cite%
{Grood-Suntay}, is used to define the rotation and translation vectors that
describe the three-dimensional patella-femoral and tibio-femoral motions.
This joint coordinate system consists of an $x$-axis that is fixed on the
femur ($\overline{i}$ is the unit vector directing the $x$-axis), a $%
z^{\prime }$-axis that is fixed on the tibia ($\overline{k}^{\prime }$ is
the unit vector directing the $z^{\prime }$-axis), and a floating axis
perpendicular to these two fixed axes ($\overline{i}\times \overline{k}%
^{\prime }$ is the unit vector directing the floating axis). The
Grood-Suntay's rotation vector includes three angular components:

\begin{enumerate}
\item $\alpha $ is the \textit{flexion-extension angle} that characterize
the rotation occured around the femoral fixed axis $O_{F}x$;

\item $\beta -\pi /2$ is the \textit{adduction-abduction angle} for the
right knee, that characterizes the rotation arround the floating axis (this
is the case studied by us). For the left knee, the adduction-abduction angle
is given by $\pi /2-\beta $. In both cases the angle $\beta $ is called the 
\textit{varus-valgus angle}.

\item $\gamma $ is the \textit{internal-external tibial angle} that
characterize the rotation which occurs about the tibial fixed axis $%
O_{T}z^{\prime }$.
\end{enumerate}

Using this joint coordinate system, the rotation vector (describing the
orientation of the tibial coordinate system with respect to the femoral
coordinate system) is defined as:

{\small 
\begin{equation*}
\overline{\theta }=-\alpha \overline{i}-\beta \left[ \overline{i}\times 
\overline{k}^{\prime }\right] -\gamma \overline{k}^{\prime }.
\end{equation*}%
}

This rotation vector can be transformed to the femoral coordinate system,
and then, it can be differentiated with respect to time to yield the angular
velocity and angular acceleration vectors of the tibia with respect to the
femur.

In this analysis, it is assumed that the femur is fixed while the tibia is
moving. The transformations of these two Grood-Suntay coordinate systems are
given by (see \cite{Grood-Suntay}, \cite{Hefzy})%
\begin{equation*}
\overline{r}=\overline{r}_{O_{T}}+\overline{r}^{\prime }=\overline{r}%
_{O_{T}}+\mathfrak{R\cdot }\overline{r}_{t=0},
\end{equation*}%
where the vector $\overline{r}_{O_{T}}$ is the position vector which locates
the origin of the tibial coordinate system with respect to the femoral
coordinate system (the tibial origin $O_{T}$ is considered in the center of
tibia), $\overline{r}_{t=0}$ describes the position vector of an arbitrary
point $P$ with respect to the tibial coordinate system (at the initial
moment $t=0$; note that the vector $\overline{r}_{t=0}$ is constant in
time), and $\overline{r}^{\prime }$ is the position vector (at an arbitrary
moment $t$) of the same point $P$ with respect to the tibial coordinate
system. As usual in the mechanics of the rigid bodies, the position vector $%
\overline{r}^{\prime }$ of the point $P$ is expressed by the rotation formula%
\begin{equation*}
\overline{r}^{\prime }=\mathfrak{R\cdot }\overline{r}_{t=0},
\end{equation*}%
where the matrix of rotation%
\begin{equation}
\mathfrak{R}{=}\left( 
\begin{array}{ccc}
R_{11} & R_{12} & R_{13} \\ 
R_{21} & R_{22} & R_{23} \\ 
R_{31} & R_{32} & R_{33}%
\end{array}%
\right)  \label{matrix of rotation R}
\end{equation}%
has the entries (see \cite{Grood-Suntay}, \cite{Hefzy})%
\begin{equation*}
R_{11}=\sin {\scriptsize \beta }\cos {\scriptsize \gamma ,}\quad R_{12}=\sin 
{\scriptsize \beta }\sin {\scriptsize \gamma ,\quad }R_{13}=\cos 
{\scriptsize \beta ,}
\end{equation*}%
\begin{equation*}
R_{21}={\scriptsize -}\cos {\scriptsize \alpha }\sin {\scriptsize \gamma -}%
\sin {\scriptsize \alpha }\cos {\scriptsize \beta }\cos {\scriptsize \gamma ,%
}
\end{equation*}%
\begin{equation*}
R_{22}=\cos {\scriptsize \alpha }\cos {\scriptsize \gamma -}\sin 
{\scriptsize \alpha }\cos {\scriptsize \beta }\sin {\scriptsize \gamma
,\quad }R_{23}=\sin {\scriptsize \alpha }\sin {\scriptsize \beta ,}
\end{equation*}%
\begin{equation*}
R_{31}=\sin {\scriptsize \alpha }\sin {\scriptsize \gamma -}\cos 
{\scriptsize \alpha }\cos {\scriptsize \beta }\cos {\scriptsize \gamma ,}
\end{equation*}%
\begin{equation*}
R_{32}={\scriptsize -}\sin {\scriptsize \alpha }\cos {\scriptsize \gamma -}%
\cos {\scriptsize \alpha }\cos {\scriptsize \beta }\sin {\scriptsize \gamma
,\quad }R_{33}=\cos {\scriptsize \alpha }\sin {\scriptsize \beta .}
\end{equation*}%
The rotation matrix $\mathfrak{R}$ is generated by multiplying the three
matrices generated by the rotations arround the axes, namely:

\begin{enumerate}
\item a rotation of angle $\alpha $ arround the axis $O_{F}x$ (from the top
of the vector $\overline{i}$ the rotation is clock-wise):%
\begin{equation*}
\mathcal{R}_{{\scriptsize O}_{F}x,\text{ }\alpha }:{=}\left( 
\begin{array}{ccc}
1 & 0 & 0 \\ 
0 & \cos \alpha & \sin \alpha \\ 
0 & -\sin \alpha & \cos \alpha%
\end{array}%
\right) ;
\end{equation*}

\item a rotation of angle $\dfrac{\pi }{2}-\beta $ arround the floating axis
(from the top of the vector $\overline{i}\times \overline{k}^{\prime }$ the
rotation is clock-wise):%
\begin{equation*}
\mathcal{R}_{\text{ {\scriptsize floating axis}},\text{ }\beta }:{=}\left( 
\begin{array}{ccc}
\sin \beta & 0 & \cos \beta \\ 
0 & 1 & 0 \\ 
-\cos \beta & 0 & \sin \beta%
\end{array}%
\right) ;
\end{equation*}

\item a rotation of angle $\gamma $ arround the axis $O_{T}z^{\prime }$
(from the top of the vector $\overline{k}^{\prime }$ the rotation is
clock-wise):%
\begin{equation*}
\mathcal{R}_{{\scriptsize O}_{T}z^{\prime },\text{ }\gamma }:{=}\left( 
\begin{array}{ccc}
\cos \gamma & \sin \gamma & 0 \\ 
-\sin \gamma & \cos \gamma & 0 \\ 
0 & 0 & 1%
\end{array}%
\right) .
\end{equation*}%
In other words, we have the matrix equality%
\begin{equation*}
\mathfrak{R}{=}\mathcal{R}_{{\scriptsize O}_{F}x,\text{ }\alpha }\cdot 
\mathcal{R}_{\text{ {\scriptsize floating axis}},\text{ }\beta }\cdot 
\mathcal{R}_{{\scriptsize O}_{T}z^{\prime },\text{ }\gamma }.
\end{equation*}
\end{enumerate}

The equations that govern the three-dimensional motion of the tibia with
respect to the femur are the second order differential Newton's and Euler's
equations of motion.

The classical Newton's equations are written in scalar form, with respect to
the femoral fixed system of axes, as:%
\begin{equation*}
\left\{ 
\begin{array}{l}
F_{x}^{\text{ext}}+G_{x}=m\ddot{x}_{O_{F}}\medskip \\ 
F_{y}^{\text{ext}}+G_{y}=m\ddot{y}_{O_{F}}\medskip \\ 
F_{z}^{\text{ext}}+G_{z}=m\ddot{z}_{O_{F}},%
\end{array}%
\right.
\end{equation*}%
where

\begin{itemize}
\item $\overline{F}^{\text{ext}}=F_{x}^{\text{ext}}\cdot \overline{i}+F_{y}^{%
\text{ext}}\cdot \overline{j}+F_{z}^{\text{ext}}\cdot \overline{k}$ is the
sum of all external forces applied to the tibia-femur (contact forces in
knee, ligamentous forces in knee etc.);

\item $\overline{G}=G_{x}\cdot \overline{i}+G_{y}\cdot \overline{j}%
+G_{z}\cdot \overline{k}$ is the weight of the leg of mass $m$; the mass of
the leg was taken in this experiment as $m=4.0$ $\unit{kg}$ (see Hefzy and
Abdel-Rahman \cite{Hefzy}).

\item $\overline{a}_{O_{F}}=\ddot{x}_{O_{F}}\cdot \overline{i}+\ddot{y}%
_{O_{F}}\cdot \overline{j}+\ddot{z}_{O_{F}}\cdot \overline{k}$ is the
acceleration of the center of mass of the leg (which is fixed in the origin
of the femoral coordinate system $O_{F}$).
\end{itemize}

From the perspective of the rigid body mechanics, the Euler's equations of
motion are written in the scalar form as:%
\begin{equation}
\left\{ 
\begin{array}{l}
I_{x^{\prime }x^{\prime }}\dot{\omega}_{x^{\prime }}+\left( I_{z^{\prime
}z^{\prime }}-I_{y^{\prime }y^{\prime }}\right) \omega _{y^{\prime }}\omega
_{z^{\prime }}=M_{x^{\prime }}^{\text{ext}}\medskip \\ 
I_{y^{\prime }y^{\prime }}\dot{\omega}_{y^{\prime }}+\left( I_{x^{\prime
}x^{\prime }}-I_{z^{\prime }z^{\prime }}\right) \omega _{z^{\prime }}\omega
_{x^{\prime }}=M_{y^{\prime }}^{\text{ext}}\medskip \\ 
I_{z^{\prime }z^{\prime }}\dot{\omega}_{z^{\prime }}+\left( I_{y^{\prime
}y^{\prime }}-I_{x^{\prime }x^{\prime }}\right) \omega _{x^{\prime }}\omega
_{y^{\prime }}=M_{z^{\prime }}^{\text{ext}},%
\end{array}%
\right.  \label{Euler equations GENERAL}
\end{equation}

where

\begin{itemize}
\item $\overline{\omega }=\omega _{x^{\prime }}\cdot \overline{i}^{\prime
}+\omega _{y^{\prime }}\cdot \overline{j}^{\prime }+\omega _{z^{\prime
}}\cdot \overline{k}^{\prime }$ is the angular velocity vector of the tibia
with respect to the femur;

\item $\overset{\bullet }{\overline{\omega }}=\dot{\omega}_{x^{\prime
}}\cdot \overline{i}^{\prime }+\dot{\omega}_{y^{\prime }}\cdot \overline{j}%
^{\prime }+\dot{\omega}_{z^{\prime }}\cdot \overline{k}^{\prime }$ is the
angular acceleration vector of the tibia with respect to the femur;

\item $\overline{M}^{\text{ext}}=M_{x^{\prime }}^{\text{ext}}\cdot \overline{%
i}^{\prime }+M_{y^{\prime }}^{\text{ext}}\cdot \overline{j}^{\prime
}+M_{z^{\prime }}^{\text{ext}}\cdot \overline{k}^{\prime }$ is the sum of
the moments (the torques) of all external forces acting on the tibia around
the $x^{\prime }$-, $y^{\prime }$-, and $z^{\prime }$-axes;

\item $I_{x^{\prime }x^{\prime }}$, $I_{y^{\prime }y^{\prime }}$ and $%
I_{z^{\prime }z^{\prime }}$ are the principal moments of inertia of the leg
about its centroidal principal system of axes. The inertial parameters were
estimated using anthropometric data as%
\begin{equation*}
I_{x^{\prime }x^{\prime }}=0.0672\unit{kg}\cdot \unit{m}^{2},\qquad
I_{y^{\prime }y^{\prime }}=0.0672\unit{kg}\cdot \unit{m}^{2},\qquad
\!\!\!\!\!I_{z^{\prime }z^{\prime }}=0.005334\unit{kg}\cdot \unit{m}^{2}.
\end{equation*}%
Note that, in this analysis, the leg was assumed to be a right cylinder. For
more details, please see Hefzy and Abdel-Rahman \cite{Hefzy}, and references
therein.
\end{itemize}

\begin{remark}
If we use the matrix equality%
\begin{equation*}
\left( 
\begin{array}{ccc}
\overline{i}^{\prime }, & \overline{j}^{\prime }, & \overline{k}^{\prime }%
\end{array}%
\right) =\left( 
\begin{array}{ccc}
\overline{i}, & \overline{j}, & \overline{k}%
\end{array}%
\right) \cdot \mathfrak{R},
\end{equation*}%
then we deduce that the rotation vector $\overline{\theta }$ takes in the
femoral system of axes the following form:%
\begin{eqnarray}
\overline{\theta } &=&\left( -\alpha -\gamma \cos \beta \right) \cdot 
\overline{i}+\left( -\beta \cos \alpha -\gamma \sin \alpha \sin \beta
\right) \cdot \overline{j}+  \label{angular rotation vector} \\
&&+\left( \beta \sin \alpha -\gamma \cos \alpha \sin \beta \right) \cdot 
\overline{k}.  \notag
\end{eqnarray}%
Therefore, by differentiating the rotation vector $\overline{\theta }$, we
find the angular velocity vector%
\begin{equation*}
\overline{\omega }=\frac{d\overline{\theta }}{dt}=\omega _{x}\cdot \overline{%
i}+\omega _{y}\cdot \overline{j}+\omega _{z}\cdot \overline{k}.
\end{equation*}%
Now, by rotating the angular velocity vector $\overline{\omega }$, via the
formula%
\begin{equation*}
\left( 
\begin{array}{ccc}
\omega _{x^{\prime }}, & \omega _{y^{\prime }}, & \omega _{z^{\prime }}%
\end{array}%
\right) =\left( 
\begin{array}{ccc}
\omega _{x}, & \omega _{y}, & \omega _{z}%
\end{array}%
\right) \cdot \mathfrak{R},
\end{equation*}%
we obtain the equalities (see also \emph{\cite{Hefzy}}):%
\begin{equation*}
\begin{array}{l}
\omega _{x^{\prime }}=-\dot{\alpha}\sin \beta \cos \gamma -\dot{\alpha}\beta
\cos \beta \cos \gamma +\dot{\alpha}\gamma \sin \beta \sin \gamma +\dot{\beta%
}\sin \gamma +\dot{\beta}\gamma \cos \gamma ,\medskip \\ 
\omega _{y^{\prime }}=-\dot{\alpha}\sin \beta \sin \gamma -\dot{\alpha}\beta
\cos \beta \sin \gamma -\dot{\alpha}\gamma \sin \beta \cos \gamma -\dot{\beta%
}\cos \gamma +\dot{\beta}\gamma \sin \gamma ,\medskip \\ 
\omega _{z^{\prime }}=-\dot{\alpha}\cos \beta +\dot{\alpha}\beta \sin \beta -%
\dot{\gamma}.%
\end{array}%
\end{equation*}
\end{remark}

\section{Approximate values of the components of the total external torque
and angular velocity}

Using the three-dimensional mathematical model for the tibio-femoral motion
due to Luh et al. (see \cite{Luh}), Apkarian et al. (see \cite{Apkarian})
practically measured in a lab frame (using some specialized devices) the
total external torques acting on the tibia of a subject (A) having the
stride length $1.41\unit{m}$ and the speed of walking $1.21\unit{m}/\unit{s}$%
. Obviously, the time of the gait cycle for the subject (A), studied in the
paper \cite{Apkarian}, is $T=1.1652\unit{s}$.

Let us consider for the subject (A) the following gait cycle intermediate
moments:%
\begin{equation*}
t_{0}=0\text{ }(\text{0\%}),\quad t_{1}=T/4=0.2913\text{ }(\text{25\%}%
),\quad t_{2}=T/2=0.5826\text{ }(\text{50\%}),
\end{equation*}%
\begin{equation*}
t_{3}=3T/4=0.8739\text{ }(\text{75\%}),\quad t_{4}=T=1.1652\text{ }(\text{%
100\%}).
\end{equation*}

Looking now at the three graphs of the knee torque components (computed in $%
\unit{N}\cdot \unit{m}$ in the paper \cite{Apkarian}, pp. 153), we see that
the total external torque vector of the knee 
\begin{eqnarray*}
\overline{M}^{\text{ext}} &=&M_{x}^{\text{ext}}\cdot \overline{i}+M_{y}^{%
\text{ext}}\cdot \overline{j}+M_{z}^{\text{ext}}\cdot \overline{k}:=\bigskip
\\
&:&=-M_{\text{valgus}}\cdot \overline{i}-M_{\text{extension}}\cdot \overline{%
j}+M_{\text{internal}}\cdot \overline{k}
\end{eqnarray*}%
has the following approximate interpolation components:%
\begin{equation*}
\begin{tabular}{||l||l||l||}
\hline\hline
$M_{x}^{\text{ext}}(t_{0})\approx 7.5,$ & $M_{y}^{\text{ext}}(t_{0})\approx
7.5,$ & $M_{z}^{\text{ext}}(t_{0})=0,$ \\ \hline\hline
$M_{x}^{\text{ext}}(t_{1})\approx -40,$ & $M_{y}^{\text{ext}}(t_{1})=0,$ & $%
M_{z}^{\text{ext}}(t_{1})\approx 5,$ \\ \hline\hline
$M_{x}^{\text{ext}}(t_{2})\approx -15,$ & $M_{y}^{\text{ext}}(t_{2})=0,$ & $%
M_{z}^{\text{ext}}(t_{2})=0,$ \\ \hline\hline
$M_{x}^{\text{ext}}(t_{3})=0,$ & $M_{y}^{\text{ext}}(t_{3})=0,$ & $M_{z}^{%
\text{ext}}(t_{3})\approx -5,$ \\ \hline\hline
$M_{x}^{\text{ext}}(t_{4})\approx 7.5,$ & $M_{y}^{\text{ext}}(t_{4})\approx
15,$ & $M_{z}^{\text{ext}}(t_{4})=0.$ \\ \hline\hline
\end{tabular}%
\end{equation*}

At the same time, looking at the approximate values of the three knee angles
appearing in the graphs from \cite{Apkarian} (pp. 150), we observe that the
rotation vector%
\begin{equation*}
\overline{\theta }=\theta _{x}\cdot \overline{i}+\theta _{y}\cdot \overline{j%
}+\theta _{z}\cdot \overline{k}:=-\theta _{\text{varus}}\cdot \overline{i}%
-\theta _{\text{flexion}}\cdot \overline{j}-\theta _{\text{external}}\cdot 
\overline{k}
\end{equation*}%
has the interpolation components (the angles are presented here in radians;
from the top of the vectors $\overline{i}$, $\overline{j}$, $\overline{k}$,
these angles are measured counter-clock-wise):%
\begin{equation*}
\begin{tabular}{||l||l||}
\hline\hline
$\theta _{x}(t_{0})=0\text{ }(=0^{\circ }),$ & $\theta _{y}(t_{0})\approx
-0.0872\text{ }(=-5^{\circ }),$ \\ \hline\hline
$\theta _{x}(t_{1})\approx -0.0872\text{ }(=-5^{\circ }),$ & $\theta
_{y}(t_{1})\approx -0.3490\text{ }(=-20^{\circ }),$ \\ \hline\hline
$\theta _{x}(t_{2})\approx -0.0872\text{ }(=-5^{\circ }),$ & $\theta
_{y}(t_{2})\approx -0.3490\text{ }(=-20^{\circ }),$ \\ \hline\hline
$\theta _{x}(t_{3})\approx -0.1745\text{ }(=-10^{\circ }),$ & $\theta
_{y}(t_{3})\approx -1.1344\text{ }(=-65^{\circ }),$ \\ \hline\hline
$\theta _{x}(t_{4})=0\text{ }(=0^{\circ }),$ & $\theta _{y}(t_{4})\approx
-0.0872\text{ }(=-5^{\circ }),$ \\ \hline\hline
\end{tabular}%
\end{equation*}%
\begin{equation*}
\begin{tabular}{||l||}
\hline\hline
$\theta _{z}(t_{0})\approx -0.1745\text{ }(=-10^{\circ }),$ \\ \hline\hline
$\theta _{z}(t_{1})\approx -0.0872\text{ }(=-5^{\circ }),$ \\ \hline\hline
$\theta _{z}(t_{2})\approx -0.1745\text{ }(=-10^{\circ }),$ \\ \hline\hline
$\theta _{z}(t_{3})\approx -0.0872\text{ }(=-5^{\circ }),$ \\ \hline\hline
$\theta _{z}(t_{4})\approx -0.1745\text{ }(=-10^{\circ }).$ \\ \hline\hline
\end{tabular}%
\end{equation*}

\begin{remark}
Solving numerically the angle system (see the rotation vector \emph{(\ref%
{angular rotation vector})})%
\begin{equation*}
\left\{ 
\begin{array}{l}
-\alpha -\gamma \cos \beta =\theta _{x}\medskip \\ 
-\beta \cos \alpha -\gamma \sin \alpha \sin \beta =\theta _{y}\medskip \\ 
\beta \sin \alpha -\gamma \cos \alpha \sin \beta =\theta _{z},%
\end{array}%
\right.
\end{equation*}%
we find the following intermediate Grood-Suntay's angles for knee:%
\begin{equation*}
\begin{tabular}{||l||l||l||}
\hline\hline
$\alpha (t_{0})\approx -0.5786,$ & $\beta (t_{0})\approx 0.1684,$ & $\gamma
(t_{0})\approx 0.5869,$ \\ \hline\hline
$\alpha (t_{1})\approx -0.0760,$ & $\beta (t_{1})\approx 0.3546,$ & $\gamma
(t_{1})\approx 0.1740,$ \\ \hline\hline
$\alpha (t_{2})\approx -0.1851,$ & $\beta (t_{2})\approx 0.3751,$ & $\gamma
(t_{2})\approx 0.2927,$ \\ \hline\hline
$\alpha (t_{3})\approx 0.0859,$ & $\beta (t_{3})\approx 1.1227,$ & $\gamma
(t_{3})\approx 0.2044,$ \\ \hline\hline
$\alpha (t_{4})\approx -0.5786,$ & $\beta (t_{4})\approx 0.1684,$ & $\gamma
(t_{4})\approx 0.5869.$ \\ \hline\hline
\end{tabular}%
\end{equation*}
\end{remark}

Now, using the Lagrange polynomial of interpolation generally described by%
\begin{eqnarray*}
P(x) &=&y_{0}\cdot \frac{x-x_{1}}{x_{0}-x_{1}}\cdot \frac{x-x_{2}}{%
x_{0}-x_{2}}\cdot \frac{x-x_{3}}{x_{0}-x_{3}}\cdot \frac{x-x_{4}}{x_{0}-x_{4}%
}+ \\
&&+y_{1}\cdot \frac{x-x_{0}}{x_{1}-x_{0}}\cdot \frac{x-x_{2}}{x_{1}-x_{2}}%
\cdot \frac{x-x_{3}}{x_{1}-x_{3}}\cdot \frac{x-x_{4}}{x_{1}-x_{4}}+ \\
&&+y_{2}\cdot \frac{x-x_{0}}{x_{2}-x_{0}}\cdot \frac{x-x_{1}}{x_{2}-x_{1}}%
\cdot \frac{x-x_{3}}{x_{2}-x_{3}}\cdot \frac{x-x_{4}}{x_{2}-x_{4}}+ \\
&&+y_{3}\cdot \frac{x-x_{0}}{x_{3}-x_{0}}\cdot \frac{x-x_{1}}{x_{3}-x_{1}}%
\cdot \frac{x-x_{2}}{x_{3}-x_{2}}\cdot \frac{x-x_{4}}{x_{3}-x_{4}}+ \\
&&+y_{4}\cdot \frac{x-x_{0}}{x_{4}-x_{0}}\cdot \frac{x-x_{1}}{x_{4}-x_{1}}%
\cdot \frac{x-x_{2}}{x_{4}-x_{2}}\cdot \frac{x-x_{3}}{x_{4}-x_{3}},
\end{eqnarray*}%
we can construct the approximate Apkarian's angle functions $\theta _{x}(t)$%
, $\theta _{y}(t)$ and $\theta _{z}(t)$. Afterward, by differentiating the
angle functions $\theta _{x}(t)$, $\theta _{y}(t)$ and $\theta _{z}(t)$, we
get the components of the angular velocity in our gait cycle moments ($\unit{%
rad}/\unit{s})$:%
\begin{equation*}
\begin{tabular}{||l||l||l||}
\hline\hline
$\omega _{x}(t_{0})\approx -1.0981,$ & $\omega _{y}(t_{0})\approx -5.6920,$
& $\omega _{z}(t_{0})\approx 1.5984,$ \\ \hline\hline
$\omega _{x}(t_{1})\approx 0.0999,$ & $\omega _{y}(t_{1})\approx 1.1983,$ & $%
\omega _{z}(t_{1})\approx -0.3997,$ \\ \hline\hline
$\omega _{x}(t_{2})\approx -0.1998,$ & $\omega _{y}(t_{2})\approx -1.7974,$
& $\omega _{z}(t_{2})\approx -0.0002,$ \\ \hline\hline
$\omega _{x}(t_{3})\approx -0.1997,$ & $\omega _{y}(t_{3})\approx -2.0970,$
& $\omega _{z}(t_{3})\approx 0.3993,$ \\ \hline\hline
$\omega _{x}(t_{4})\approx 1.8975,$ & $\omega _{y}(t_{4})\approx 12.8820,$ & 
$\omega _{z}(t_{4})\approx -1.5983.$ \\ \hline\hline
\end{tabular}%
\end{equation*}%
Consequently, via the rotation formulas%
\begin{equation*}
\begin{array}{l}
\text{X}^{^{\prime }}=R_{11}\text{X}+R_{21}\text{Y}+R_{31}\text{Z},\medskip
\\ 
\text{Y}^{^{\prime }}=R_{12}\text{X}+R_{22}\text{Y}+R_{32}\text{Z},\medskip
\\ 
\text{Z}^{^{\prime }}=R_{13}\text{X}+R_{23}\text{Y}+R_{33}\text{Z},%
\end{array}%
\end{equation*}%
where $R_{ij}$ are given by the expressions from (\ref{matrix of rotation R}%
), we find the following intermediate vector components (with respect to
tibia coordinate system):\medskip

$\blacksquare $ $-$ components of the total external torque:%
\begin{equation*}
\begin{tabular}{||l||l||l||}
\hline\hline
$M_{x^{\prime }}^{\text{ext}}(t_{0})\approx 0.9361,$ & $M_{y^{\prime }}^{%
\text{ext}}(t_{0})\approx 8.1637,$ & $M_{z^{\prime }}^{\text{ext}%
}(t_{0})\approx 6.7065,$ \\ \hline\hline
$M_{x^{\prime }}^{\text{ext}}(t_{1})\approx -18.3490,$ & $M_{y^{\prime }}^{%
\text{ext}}(t_{1})\approx -2.8400,$ & $M_{z^{\prime }}^{\text{ext}%
}(t_{1})\approx -35.7800,$ \\ \hline\hline
$M_{x^{\prime }}^{\text{ext}}(t_{2})\approx -5.2618,$ & $M_{y^{\prime }}^{%
\text{ext}}(t_{2})\approx -1.5857,$ & $M_{z^{\prime }}^{\text{ext}%
}(t_{2})\approx -13.9570,$ \\ \hline\hline
$M_{x^{\prime }}^{\text{ext}}(t_{3})\approx 2.0263,$ & $M_{y^{\prime }}^{%
\text{ext}}(t_{3})\approx 0.8581,$ & $M_{z^{\prime }}^{\text{ext}%
}(t_{3})\approx -4.4898,$ \\ \hline\hline
$M_{x^{\prime }}^{\text{ext}}(t_{4})\approx 0.8255,$ & $M_{y^{\prime }}^{%
\text{ext}}(t_{4})\approx 15.6310,$ & $M_{z^{\prime }}^{\text{ext}%
}(t_{4})\approx 6.0191.$ \\ \hline\hline
\end{tabular}%
\end{equation*}

$\blacksquare $ $-$ components of the angular velocity:%
\begin{equation*}
\begin{tabular}{||l||l||l||}
\hline\hline
$\omega _{x^{\prime }}(t_{0})\approx -1.6519,$ & $\omega _{y^{\prime
}}(t_{0})\approx -5.7721,$ & $\omega _{z^{\prime }}(t_{0})\approx -0.3366,$
\\ \hline\hline
$\omega _{x^{\prime }}(t_{1})\approx 0.2847,$ & $\omega _{y^{\prime
}}(t_{1})\approx 1.2324,$ & $\omega _{z^{\prime }}(t_{1})\approx -0.0762,$
\\ \hline\hline
$\omega _{x^{\prime }}(t_{2})\approx 0.1451,$ & $\omega _{y}(t_{2})\approx
-1.8010,$ & $\omega _{z^{\prime }}(t_{2})\approx -0.0647,$ \\ \hline\hline
$\omega _{x^{\prime }}(t_{3})\approx 0.1623,$ & $\omega _{y^{\prime
}}(t_{3})\approx -2.1350,$ & $\omega _{z^{\prime }}(t_{3})\approx 0.1098,$
\\ \hline\hline
$\omega _{x^{\prime }}(t_{4})\approx 1.6574,$ & $\omega _{y^{\prime
}}(t_{4})\approx 13.0050,$ & $\omega _{z^{\prime }}(t_{4})\approx 0.4656.$
\\ \hline\hline
\end{tabular}%
\end{equation*}

Now, using the preceding Tables and the well-known method of multiple
regression from Statistics, by numerical computations, we get the following
linear approximations for the torque components $M_{x^{\prime }}^{\text{ext}%
} $, $M_{y^{\prime }}^{\text{ext}}$ and $M_{z^{\prime }}^{\text{ext}}$, with
respect to the angular velocity components $\omega _{x^{\prime }}$, $\omega
_{y^{\prime }}$ and $\omega _{z^{\prime }}$ (via the classical method of
least squares):%
\begin{equation}
\begin{array}{l}
M_{x^{\prime }}^{\text{ext}}\approx -16.9430\cdot \omega _{x^{\prime
}}-0.5003\cdot \omega _{y^{\prime }}+80.9290\cdot \omega _{z^{\prime
}}-3.0709,\medskip \\ 
M_{y^{\prime }}^{\text{ext}}\approx -15.1720\cdot \omega _{x^{\prime
}}+1.5740\cdot \omega _{y^{\prime }}+34.8270\cdot \omega _{z^{\prime
}}+3.7511,\medskip \\ 
M_{z^{\prime }}^{\text{ext}}\approx -39.7610\cdot \omega _{x^{\prime
}}+0.9071\cdot \omega _{y^{\prime }}+140.8000\cdot \omega _{z^{\prime
}}-7.1266.%
\end{array}
\label{approximative-torques}
\end{equation}

\section{From Euler's equations of motion of knee to jet "electromagnetic"
Yang-Mills energy}

Note that, in our jet geometrical approach, we can regard the Euler's
equations of motion (\ref{Euler equations GENERAL}) as a dynamical system on
the $1$-jet space $J^{1}(\mathbb{R},\mathbb{R}^{3})$. The coordinates on the 
$1$-jet space $J^{1}(\mathbb{R},\mathbb{R}^{3})$ are considered as being%
\begin{equation*}
(t,\omega _{x^{\prime }},\omega _{y^{\prime }},\omega _{z^{\prime }},\dot{%
\omega}_{x^{\prime }},\dot{\omega}_{y^{\prime }},\dot{\omega}_{z^{\prime }}).
\end{equation*}

In such a context, taking the particular case from the Hefzy --
Abdel-Rahman's paper \cite{Hefzy}, when we have%
\begin{equation*}
I_{x^{\prime }x^{\prime }}=I_{y^{\prime }y^{\prime }}=0.0672,\quad
I_{z^{\prime }z^{\prime }}=0.0053,
\end{equation*}%
and setting the components of the total external torque by the expressions (%
\ref{approximative-torques}), then the Euler's equations of motion (\ref%
{Euler equations GENERAL}) for the subject (A) can be rewritten in the
following particular form:%
\begin{equation}
\left\{ 
\begin{array}{lll}
\dot{\omega}_{x^{\prime }} & = & 0.9211\cdot \omega _{y^{\prime }}\cdot
\omega _{z^{\prime }}-252.1279\cdot \omega _{x^{\prime }}-7.4449\cdot \omega
_{y^{\prime }}+\medskip  \\ 
&  & +1204.3005\cdot \omega _{z^{\prime }}-45.6979\medskip  \\ 
\dot{\omega}_{y^{\prime }} & = & -0.9211\cdot \omega _{x^{\prime }}\cdot
\omega _{z^{\prime }}-225.7738\cdot \omega _{x^{\prime }}+23.4226\cdot
\omega _{y^{\prime }}+\medskip  \\ 
&  & +518.2589\cdot \omega _{z^{\prime }}+55.8199\medskip  \\ 
\dot{\omega}_{z^{\prime }} & = & -7502.0754\cdot \omega _{x^{\prime
}}+171.1509\cdot \omega _{y^{\prime }}+\medskip  \\ 
&  & +26566.0377\cdot \omega _{z^{\prime }}-1344.6415.%
\end{array}%
\right.   \label{Euler PARTICULAR}
\end{equation}

Consequently, via the Theorem \ref{MainTh} applied to the first order ODE
system (\ref{Euler PARTICULAR}), we assert that the Lagrangian geometrical
behavior on the $1$-jet space $J^{1}(\mathbb{R},\mathbb{R}^{3})$ of the
above knee dynamical system can be described by the following result:

\begin{corollary}
\begin{enumerate}
\item[\emph{(1)}] The \textbf{canonical nonlinear connection} on $J^{1}(%
\mathbb{R},\mathbb{R}^{3})$\ \textbf{produced by the jet dynamical system of
knee }\emph{(\ref{Euler PARTICULAR})} has the local components 
\begin{equation*}
\mathring{\Gamma}=\left( \boldsymbol{\mathring{M}}_{(1)1}^{(i)}=0,%
\boldsymbol{\mathring{N}}_{(1)j}^{(i)}\right) ,
\end{equation*}%
where the matrix of the spatial nonlinear connection%
\begin{equation*}
\mathbf{\mathring{N}}_{(1)}=\left( 
\begin{array}{ccc}
\boldsymbol{\mathring{N}}_{(1)1}^{(1)} & \boldsymbol{\mathring{N}}%
_{(1)2}^{(1)} & \boldsymbol{\mathring{N}}_{(1)3}^{(1)}\medskip \\ 
\boldsymbol{\mathring{N}}_{(1)1}^{(2)} & \boldsymbol{\mathring{N}}%
_{(1)2}^{(2)} & \boldsymbol{\mathring{N}}_{(1)3}^{(2)}\medskip \\ 
\boldsymbol{\mathring{N}}_{(1)1}^{(3)} & \boldsymbol{\mathring{N}}%
_{(1)2}^{(3)} & \boldsymbol{\mathring{N}}_{(1)3}^{(3)}%
\end{array}%
\right) =\left( 
\begin{array}{ccc}
0 & -\boldsymbol{\mathring{F}}_{(1)2}^{(1)} & -\boldsymbol{\mathring{F}}%
_{(1)3}^{(1)}\medskip \\ 
\boldsymbol{\mathring{F}}_{(1)2}^{(1)} & 0 & -\boldsymbol{\mathring{F}}%
_{(2)3}^{(1)}\medskip \\ 
\boldsymbol{\mathring{F}}_{(1)3}^{(1)} & \boldsymbol{\mathring{F}}%
_{(2)3}^{(1)} & 0%
\end{array}%
\right)
\end{equation*}%
has the entries%
\begin{equation*}
\boldsymbol{\mathring{F}}_{(1)2}^{(1)}=0.9211\cdot \omega _{z^{\prime
}}+109.1644,\quad \boldsymbol{\mathring{F}}_{(1)3}^{(1)}=0.4605\cdot \omega
_{y^{\prime }}+4353.1879,
\end{equation*}%
\begin{equation*}
\boldsymbol{\mathring{F}}_{(2)3}^{(1)}=-0.4605\cdot \omega _{x^{\prime
}}+173.5540.
\end{equation*}

\item[\emph{(2)}] All adapted components of the \textbf{canonical Cartan
linear connection} $\mathtt{C}\mathring{\Gamma}$\ \textbf{produced by the
jet dynamical system of knee} \emph{(\ref{Euler PARTICULAR})} are zero.

\item[\emph{(3)}] The \textbf{torsion tensor} $\mathbf{\mathring{T}}$ of the
canonical Cartan linear connection $\mathtt{C}\mathring{\Gamma}$ \textbf{%
produced by the jet dynamical system of knee} \emph{(\ref{Euler PARTICULAR})}
has as adapted components the entries of the following torsion matrices:%
\begin{equation*}
\mathbf{\mathring{T}}_{(1)1}=\dfrac{\partial \mathbf{\mathring{N}}_{(1)}}{%
\partial \omega _{x^{\prime }}}=\left( 
\begin{array}{ccc}
0 & 0 & 0\medskip \\ 
0 & 0 & 0.4605\medskip \\ 
0 & -0.4605 & 0%
\end{array}%
\right) ,
\end{equation*}%
\begin{equation*}
\mathbf{\mathring{T}}_{(1)2}=\dfrac{\partial \mathbf{\mathring{N}}_{(1)}}{%
\partial \omega _{y^{\prime }}}=\left( 
\begin{array}{ccc}
0 & 0 & -0.4605\medskip \\ 
0 & 0 & 0\medskip \\ 
0.4605 & 0 & 0%
\end{array}%
\right) ,
\end{equation*}%
\begin{equation*}
\mathbf{\mathring{T}}_{(1)3}=\dfrac{\partial \mathbf{\mathring{N}}_{(1)}}{%
\partial \omega _{z^{\prime }}}=\left( 
\begin{array}{ccc}
0 & -0.9211 & 0\medskip \\ 
0.9211 & 0 & 0\medskip \\ 
0 & 0 & 0%
\end{array}%
\right) .
\end{equation*}

\item[\emph{(4)}] All adapted components of the \textbf{curvature tensor} $%
\mathbf{\mathring{R}}$ of the canonical Cartan linear connection $\mathtt{C}%
\mathring{\Gamma}$ \textbf{produced by the jet dynamical system of knee} 
\emph{(\ref{Euler PARTICULAR})} cancel.

\item[\emph{(5)}] The \textbf{geometric "electromagnetic" adapted components
produced by the jet dynamical system of knee} \emph{(\ref{Euler PARTICULAR})}
are the entries of the matrix%
\begin{equation*}
\mathbf{\mathring{F}}^{(1)}=-\mathbf{\mathring{N}}_{(1)}=\left( 
\begin{array}{ccc}
0 & \boldsymbol{\mathring{F}}_{(1)2}^{(1)} & \boldsymbol{\mathring{F}}%
_{(1)3}^{(1)}\medskip \\ 
-\boldsymbol{\mathring{F}}_{(1)2}^{(1)} & 0 & \boldsymbol{\mathring{F}}%
_{(2)3}^{(1)}\medskip \\ 
-\boldsymbol{\mathring{F}}_{(1)3}^{(1)} & -\boldsymbol{\mathring{F}}%
_{(2)3}^{(1)} & 0%
\end{array}%
\right) .
\end{equation*}

\item[\emph{(6)}] The \textbf{jet geometric "electromagnetic" Yang-Mills
energy produced by the jet dynamical system of knee} \emph{(\ref{Euler
PARTICULAR})} is given by the formula 
\begin{equation*}
\mathbf{EYM}^{\text{\emph{knee}}}(\omega _{x^{\prime }},\omega _{y^{\prime
}},\omega _{z^{\prime }})=\left[ \boldsymbol{\mathring{F}}_{(1)2}^{(1)}%
\right] ^{2}+\left[ \boldsymbol{\mathring{F}}_{(1)3}^{(1)}\right] ^{2}+\left[
\boldsymbol{\mathring{F}}_{(2)3}^{(1)}\right] ^{2}.
\end{equation*}
\end{enumerate}
\end{corollary}

\begin{proof}
The Euler's equations of motion (\ref{Euler PARTICULAR}) represent a
particular case of the jet first order ODE system (\ref{DEs1}) for $n=3$ and 
$\mathfrak{X}=\left( \boldsymbol{\mathring{X}}_{(1)}^{(i)}(\omega
_{x^{\prime }},\omega _{y^{\prime }},\omega _{z^{\prime }})\right) _{i=%
\overline{1,3}},$ where%
\begin{equation*}
\begin{array}{lll}
\boldsymbol{\mathring{X}}_{(1)}^{(1)}(\omega _{x^{\prime }},\omega
_{y^{\prime }},\omega _{z^{\prime }}) & = & 0.9211\cdot \omega _{y^{\prime
}}\cdot \omega _{z^{\prime }}-252.1279\cdot \omega _{x^{\prime
}}-7.4449\cdot \omega _{y^{\prime }}+\medskip \\ 
&  & +1204.3005\cdot \omega _{z^{\prime }}-45.6979,%
\end{array}%
\end{equation*}%
\begin{equation*}
\begin{array}{lll}
\boldsymbol{\mathring{X}}_{(1)}^{(2)}(\omega _{x^{\prime }},\omega
_{y^{\prime }},\omega _{z^{\prime }}) & = & -0.9211\cdot \omega _{x^{\prime
}}\cdot \omega _{z^{\prime }}-225.7738\cdot \omega _{x^{\prime
}}+23.4226\cdot \omega _{y^{\prime }}+\medskip \\ 
&  & +518.2589\cdot \omega _{z^{\prime }}+55.8199,%
\end{array}%
\end{equation*}%
\begin{equation*}
\begin{array}{lll}
\boldsymbol{\mathring{X}}_{(1)}^{(3)}(\omega _{x^{\prime }},\omega
_{y^{\prime }},\omega _{z^{\prime }}) & = & -7502.0754\cdot \omega
_{x^{\prime }}+171.1509\cdot \omega _{y^{\prime }}+\medskip \\ 
&  & +26566.0377\cdot \omega _{z^{\prime }}-1344.6415.%
\end{array}%
\end{equation*}%
Consequently, using the formulas%
\begin{equation*}
\begin{array}{cc}
\boldsymbol{\mathring{F}}_{(1)2}^{(1)}=\dfrac{1}{2}\left[ \dfrac{\partial 
\boldsymbol{\mathring{X}}_{(1)}^{(1)}}{\partial \omega _{y^{\prime }}}-%
\dfrac{\partial \boldsymbol{\mathring{X}}_{(1)}^{(2)}}{\partial \omega
_{x^{\prime }}}\right] , & \boldsymbol{\mathring{F}}_{(1)3}^{(1)}=\dfrac{1}{2%
}\left[ \dfrac{\partial \boldsymbol{\mathring{X}}_{(1)}^{(1)}}{\partial
\omega _{z^{\prime }}}-\dfrac{\partial \boldsymbol{\mathring{X}}_{(1)}^{(3)}%
}{\partial \omega _{x^{\prime }}}\right] ,%
\end{array}%
\end{equation*}%
\begin{equation*}
\begin{array}{cc}
\boldsymbol{\mathring{F}}_{(2)3}^{(1)}=\dfrac{1}{2}\left[ \dfrac{\partial 
\boldsymbol{\mathring{X}}_{(1)}^{(2)}}{\partial \omega _{z^{\prime }}}-%
\dfrac{\partial \boldsymbol{\mathring{X}}_{(1)}^{(3)}}{\partial \omega
_{y^{\prime }}}\right] , & \mathbf{\mathring{F}}^{(1)}=-\mathbf{\mathring{N}}%
_{(1)},%
\end{array}%
\end{equation*}%
some direct computations lead us to what we were looking for.
\end{proof}

\section{Yang-Mills energetic surfaces of constant level produced by the jet
dynamical system of knee}

Using the notations $\omega _{x^{\prime }}:=X$, $\omega _{y^{\prime }}:=Y$, $%
\omega _{z^{\prime }}:=Z$, then the jet geometric "electromagnetic"
Yang-Mills energy produced by the dynamical system of knee (\ref{Euler
PARTICULAR}) (associated to the subject (A)) takes the form 
\begin{eqnarray*}
\mathbf{EYM}^{\text{knee}}(X,Y,Z) &=&\left[ 0.9211Z+109.1644\right] ^{2}+%
\left[ 0.4605Y+4353.1879\right] ^{2}+ \\
&&+\left[ -0.4605X+173.5540\right] ^{2}.
\end{eqnarray*}

Now, let us consider the \textit{jet "electromagnetic" Yang-Mills energetic
surfaces of constant level}, which are produced by the dynamical system of
the knee of the subject (A). These energetic surfaces are defined by the
implicit equations%
\begin{equation*}
\Sigma _{k}^{\text{knee}}:\mathbf{EYM}^{\text{knee}}(X,Y,Z)=k,
\end{equation*}%
where $k\geq 0$ is a given constant real number. By direct numerical
computations, we obtain the following approximate implicit equations:%
\begin{equation*}
\Sigma _{k}^{\text{knee}}:0.2120\left( X^{\prime }\right) ^{2}+0.2120\left(
Y^{\prime }\right) ^{2}+0.8484\left( Z^{\prime }\right) ^{2}=k\Leftrightarrow
\end{equation*}%
\begin{equation*}
\Sigma _{k}^{\text{knee}}:\frac{\left( X^{\prime }\right) ^{2}}{4.7169k}+%
\frac{\left( Y^{\prime }\right) ^{2}}{4.7169k}+\frac{\left( Z^{\prime
}\right) ^{2}}{1.1786k}-1=0,\quad k\neq 0,
\end{equation*}%
where we have done the spatial translation%
\begin{equation*}
X^{\prime }=X-376.8816,\quad Y^{\prime }=Y+9453.1767,\quad Z^{\prime
}=Z+118.5152.
\end{equation*}

In conclusion, it follows that the jet "electromagnetic" Yang-Mills
energetic surface of constant level $k\geq 0$ is:

\begin{enumerate}
\item for $k=0$, it is the \textit{point} $\mathcal{C}%
(376.8816,-9453.1767,-118.5152).$

\item for $k>0$, it is an \textit{ellipsoid} centered in the point $\mathcal{%
C}$, whose semi-axes are parallel with the axes $OX,$ $OY$ and $OZ$, and
they have the lengths%
\begin{equation*}
a=b=\sqrt{4.7169k},\qquad c=\sqrt{1.1786k}.
\end{equation*}
\end{enumerate}

\textbf{Open problem. }Find the biomechanical interpretations for the
ellipsoid shapes (oblate spheroids) of the jet Yang-Mills energetic surfaces
of constant level produced by the dynamical system of knee (\ref{Euler
PARTICULAR}).

Mircea NEAGU

University Transilvania of Bra\c{s}ov,

Department of Mathematics - Informatics,

Blvd. Iuliu Maniu, no. 50, Bra\c{s}ov 500091, Romania.

\textit{E-mail:} mircea.neagu@unitbv.ro\bigskip

Mihaela Maria MARIN (student in Medical Engineering)

University Transilvania of Bra\c{s}ov, Romania.

\textit{E-mail:} marin\_mihaela\_maria@yahoo.com

\end{document}